\documentclass[11pt]{article}
\usepackage{amsmath,amssymb,amsthm} 
\usepackage{graphicx}
\usepackage{amsfonts}
\usepackage{bbold}

\def \R{\mathbf{R}}
\def \N{\mathbf{N}}

\newtheorem{lemma}{Lemma}[section]
\newtheorem{proposition}[lemma]{Proposition}
\newtheorem{theorem}[lemma]{Theorem}
\newtheorem{remark}[lemma]{Remark}
\newtheorem{Aassumption}{Assumption}
\newtheorem{definition}[lemma]{Definition}

\def\reff#1{{\rm(\ref{#1})}}
\newtheorem{Sassumption}{Assumption}

\def\reff#1{{\rm(\ref{#1})}}

\def \proof{{\noindent \bf Proof. }}

 \def\reff#1{{\rm(\ref{#1})}}

\usepackage[margin=1in]{geometry}

\begin{document}

\title{Some non monotone schemes for Hamilton-Jacobi-Bellman equations}
\author{Xavier Warin\footnote{EDF R\&D \& FiME, Laboratoire de Finance des March\'es de l'Energie} \thanks{xavier.warin@edf.fr}}

\maketitle

\begin{abstract}
We extend the theory of Barles Jakobsen \cite{barles2005error} for a class of  almost monotone schemes to solve stationary Hamilton Jacobi Bellman equations.
We show that the monotonicity of the schemes can be relaxed still leading to the convergence to the viscosity solution of the equation even if the discrete problem can only be solved with some error.
We give some examples of such numerical schemes  and show that the bounds obtained by the framework developed are not tight. At last we test the schemes.
\end{abstract}

\noindent {\bf Keywords:}  Hamilton-Jacobi-Bellman equations, stochastic control, numerical methods. 
\vspace{2mm}

\noindent {\bf MSC2010:}  49L20, 65N12.

\section*{Introduction}

We are interested in the following HJB equation arising in infinite horizon, discounted, stochastic control problems
\begin{eqnarray}
\label{eqHJB}
F(x,u, \mathbf{D}u,  \mathcal{D}^2u) = 0  \mbox{ in } \R^N,
\end{eqnarray}
with 
\begin{eqnarray}
F(x,t,p,X) & = &  \sup_{\alpha \in \mathcal{A}} \mathcal{L}^{\alpha}(x,t,p,X),  \nonumber \\
\mathcal{L}^{\alpha}(x,t,p,X) & = & -tr[ a^{\alpha}(x) X] - b^{\alpha}(x)p + c^{\alpha}(x) t - f^{\alpha}(x). \label{Loperator}
\end{eqnarray}
where $a,b,c,f$ are at least continuous functions on $\R^N \times \mathcal{A}$ with values in $S(N)$ the space of symmetric $N \times N$ matrices,
$\R^N$, $\R$ and $\R$  respectively. The space of controls $ \mathcal{A}$ is supposed to be a compact metric space.\\
Supposing  $h$ is an approximating parameter, we consider an approximation $S$ of $F$ such that the approximate function $u_h$ satisfies:
\begin{eqnarray}
\label{eqHJBMono}
S(h, x, u_h(x), [u_h]_x) = 0, x \in \R^N,
\end{eqnarray}
where  $S(h, x , r, [t]_x)$ is defined for  $x \in \R^N$, $r \in \R $, 
$t$ a function defined on $\R^N$, and $[t]_x$  is a function defined at $x$ from $t$. 
This notation was introduced by \cite{barles1991convergence} to prove that a  scheme $S$  which is non decreasing in $r$ and non increasing in $[t]_x$ is  monotone. 
When the scheme $S$ is a monotone, uniformly continuous  and a consistent approximation of  $F$  and when a discrete bounded solution $u_h$ can be found for \reff{eqHJBMono},
then $u_h$ converges to the viscosity solution
of the problem \reff{eqHJB} \cite{barles2002convergence,barles2005error}.
It is also possible to modify the consistency version proposed by \cite{barles1991convergence} and still get a scheme converging to the viscosity solution for example for some finite element schemes defined in \cite{jensen2013convergence}.\\
When the scheme is non monotone, some theory exists to get convergent schemes in Sobolev space \cite{smears2016discontinuous} but breaking monotonicity can lead to non converging schemes as shown in \cite{oberman2006convergent}, and \cite{pooley2003numerical} gave some examples of non monotone schemes
converging toward  a false solution.\\
Methods to solve HJB equations includes Finite Difference methods and Semi Lagrangian methods.
Classical Finite Difference method often can be interpreted as a Markov Chain \cite{kushner2013numerical} leading to monotone schemes.
When  $a^\alpha$ is not diagonally dominant, the requirement about monotonicity of the scheme leads to Finite Difference scheme such as in \cite{bonnans2004fast} using ideas independently developed in \cite{krylov2005rate}. As an alternative to Finite Difference schemes, Semi Lagrangian schemes of low order based on the original work in \cite{menaldi1989some} have been developed in \cite{camilli1995approximation,barles2002convergence}.\\
Monotonicity of the schemes is desirable because a convenient framework is available. Nevertheless, it is  to notice that some high order non monotone schemes have been developed and proved convergent for some first order Hamilton-Jacobi equations  for examples in \cite{augoula2000high,lions1995convergence}  or in \cite{lepsky2000spectral} using
spectral methods on a periodic domain.
In the case of first order Hamilton-Jacobi-Bellman equation, some non monotone explicit schemes were developed because of the inefficiency of monotone scheme for discontinuous initial  data and proved to be convergent in \cite{bokanowski2010convergence}. 

In this paper we  will relax too the constraint on the monotonicity of the scheme such that it can converge to the right solution.
We will suppose that this scheme is a perturbation of a monotone scheme $\hat S$ and we have in mind the schemes based on interpolation method (Semi Lagrangian scheme or Finite Difference scheme with carefully chosen directions).
This idea is not new : \cite{barles1991convergence} already have emphasized the fact that monotonicity could be relaxed. 
As for Finite Difference schemes, for first order Hamilton-Jacobi, following the ideas in \cite{abgrall2009construction,froese2013convergent}, some potentially high order scheme dubbed filtered schemes have been developed  in \cite{oberman2015filtered,bokanowski2016efficient} by blending two schemes: one of high order potentially instable and one monotone of low order that will be used near singularity of the solution.  Some  filtered schemes have been also used for second order HJB equations for example in \cite{bokanowski2016high}.
As for Semi Lagrangian, results for high order schemes interpolators are given by the Italian School  for first order Hamilton-Jacobi in \cite{falcone2013semi}. The interest in nearly monotone scheme is driven by the fact that  monotone scheme such as Semi-Lagrangian schemes with linear interpolators are converging numerically very slowly.
 Then it seems natural to try to use schemes that are potentially
fast convergent  (with  potentially  a high order of consistency) while being sure that they  converge towards the true solution even if we cannot prove a higher rate of convergence than in the monotone case.
Our goal is to develop a framework  that  could be used to develop new schemes and easily prove that they  are convergent.
 \\
Specifically we treat the second order Hamilton Jacobi Bellman stationary problem with nearly monotone schemes.
Because in some case the resolution of \reff{eqHJBMono} can  be impossible we will try to relax the equality requirement  and only assume  we can find a function $u_h$  such that on a given grid $X$ that may depend on $h$ 
\begin{eqnarray}
\label{eqHJBDis}
|S(h, x, u_h(x), [u_h]_x)| \le  \epsilon(h), x \in X,
\end{eqnarray}
where $\epsilon(h)$ is a continuous function of $h$ with $\epsilon(0) =0$.\\

If existing, the solution of such a scheme is not unique, so we consider a  constructed sequence of solution  $u_h$ of \reff{eqHJBDis} and  get bounds   proving the convergence of $u_h$ towards the right solution.

This work is motivated by the previous work by the author \cite{warin2016some} : some almost monotone Semi Lagrangian schemes with high order interpolation and truncation were developed in the time dependent case.
It was proven that such schemes could be used to estimate the solution of time dependent Hamilton-Jacobi-Bellman equations. A numerical study was achieved comparing different types of Lagrange interpolators, Bernstein approximations, and Cubic spline interpolators. The conclusion was that Gauss Lobatto Legendre (GLL) interpolators were the most interesting given the accuracy obtained compared to the CPU time used.\\
In the stationary case classical techniques involving fixed point iteration scheme cannot  be used to proved existence of a solution of the discretized problem: it leads to the idea of relaxation given by equation \eqref{eqHJBDis} and to the development of a general framework.\\
In a second part of the article,  we detail the semi-lagrangian schemes proposed by \cite{warin2016some} 
and show that they can be cast into this framework so that  some convergence properties can be derived.
In the present  article, we only focus on  GLL interpolators to study the stationary  Hamilton-Jacobi-Bellman equations in our framework.
Using a direct estimation, we besides prove that the result obtained  by the framework is not optimal. In fact, with the direct estimation, we get back the  convergence result in $O(h^{\frac{1}{4}})$ previously obtained in \cite{camilli1995approximation,barles2002convergence}.

We then develop a Finite Difference approach with interpolation. In order to use the developed framework, we have to suppose that the diffusion coefficient is 
independent of the space.
Once again the method can be cast in the framework developed and  the convergence rate obtained is not optimal. The best rate we found  is in $O(h^{\frac{1}{2}})$ which is the rate found in \cite{barles2002convergence} with the same assumptions.  As for the general case where the diffusion depends on the space, the same rate of convergence is  reached  using the Bonnans and al. \cite{bonnans2004fast} or Krylov \cite{krylov2005rate} discretization  as shown in \cite{krylov2005rate}. 
\\

In the sequel the constant $C$ may vary between lines.

\section{Main result}
We define the norm denoted $| |$ as follows: for any integer $m \ge 1$ and $z = (z_i)_i \in \R^m$, we
set  $|z|^2 = \sum_{i=1}^m z_i^2$. For a matrix $M \in \R^{n_1 \times \ n_2}$, $|M|^2 = tr[ M^t M]$ with $M^t$ the transpose of $M$.\\
If $f : \R^N \longrightarrow \R^M$  we define the semi-norms:
\begin{eqnarray}
|f|_0 = \sup_{x \in \R^N} |f(x)|, \quad  [f]_{1} = \sup_{ \begin{array}{c} x, y \in \R^N \\ x \ne y \end{array}} \frac{|f(x) -f(y)|}{|x-y|},  \nonumber
\end{eqnarray}
and
\begin{eqnarray}
|f|_1 = |f|_0 + [f]_1. \nonumber
\end{eqnarray}
$C^{0,1}(\R^N)$ stands for the set of functions $f : \R^n \longrightarrow \R$ with finite norm $|f|_{1}$, $C_b(\R^N)$ the set with
finite norm $|f|_0$.
In the sequel we make the following classical assumptions
\begin{Aassumption}
\label{refAss1}
For any $\alpha \in \mathcal{A}, a^{\alpha} = \frac{1}{2} \sigma^{\alpha} \sigma^{\alpha t}$ for some $N \times P$ matrix $\sigma^{\alpha}$.
Furthermore, there exists $\lambda$, $K$ independent of $\alpha$ such that:
\begin{eqnarray}
\label{refCoef}
 c^{\alpha} \ge \lambda  > 0, \mbox{ and }  |\sigma^{\alpha}|_1 + |b^{\alpha}|_1 + |f^{\alpha}|_1 \le K.
\end{eqnarray}
\end{Aassumption}
\begin{Aassumption}
\label{refAss2} 
The constant $\lambda$ in \ref{refAss1} satisfies $\lambda  > \sup_{\alpha}{ \frac{1}{2}  [ \sigma^{\alpha}]^2_1 + [b^{\alpha}]_1}$.
\end{Aassumption}
We just recall the well-posedness and regularity result given in \cite{barles2002convergence} with demonstrations in the references therein.
\begin{proposition}
Assume \ref{refAss1}: There exists a unique viscosity $u$ solution of \reff{eqHJB} in $C_b(\R^N)$. If $w_1$ and $w_2$  are in $C_b(\R^N)$ and are sub- and supersolution of  \reff{eqHJB} respectively, then $w_1 \le w_2$ in $\R^N$.\\
Assume \ref{refAss1}, \ref{refAss2}: There exists a unique bounded viscosity $u$ solution of \reff{eqHJB} in $C^{0,1}(\R^N)$.
\end{proposition}
\begin{remark}
Assumptions \ref{refAss1}  can be given for more general H\"older spaces, and regularity of the solution is then given in \cite{barles2005error}.
\end{remark}
Here we add some new definitions that will be helpful in the sequel. 
First we introduce the notion of $\epsilon$ monotone scheme  stating that the scheme $S$ is ``nearly'' monotone
\begin{definition}
{\bf An $\epsilon(p,K )$  monotone scheme $S$} is a scheme such that there exists  $\bar \lambda$   satisfying:
\begin{itemize}
\item for every $h > 0$, $x\in \R^N$, $r \in \R$, for every function $w\in C^{0,1}(\R^N)$ , $v \in C_b(\R^N)$  such that $ v \ge w$:
\begin{eqnarray}
\label{monoton1}
S(h, x,r , [v]_x) \le   S(h, x, r, [w]_x) +  K |w|_1 h^p,
\end{eqnarray}
\item for every $h > 0$, $x\in \R^N$, $r \in \R$, for every function $w\in C_b(\R^N)$ , $v \in C^{0,1}(\R^N)$  such that $ v \ge w$:
\begin{eqnarray}
\label{monoton1bis}
S(h, x,r , [v]_x) \le   S(h, x, r, [w]_x) +  K |v|_1 h^p,
\end{eqnarray}
\item for every $h > 0$, $x\in \R^N$, $r \in \R$, $m \ge 0$,  for every function  $u \in C_b(\R^N)$:
\begin{eqnarray}
\label{monoton2}
S(h, x, r+m, [u+m]_x)  \ge  S(h, x, r, [u]_x) + \bar \lambda m.
\end{eqnarray}
\end{itemize}
\end{definition}
Because it is sometimes difficult or impossible to prove that a discrete scheme has a solution, we relax the notion of solution as we did in \reff{eqHJBDis}:
\begin{definition}
  Relaxing the notion of solution, subsolution and supersolution we introduce:
 \begin{itemize}
    \item
{\bf An $\epsilon(c)$ solution $u$ of scheme $S$} is a continuous function which satisfies:
\begin{eqnarray}
|S(h, x, u(x), [u]_x)| < c. \nonumber
\end{eqnarray} 
\item
{\bf An $\epsilon(c)$ subsolution (supersolution) $u$ of scheme $S$} is a continuous function which satisfies:
\begin{eqnarray}
S(h, x, u(x), [u]_x) < c   \quad \quad ( S(h, x, u(x), [u]_x) > - c ). \nonumber
\end{eqnarray}
\end{itemize}
\end{definition}

On the scheme we make the following assumptions:
\begin{Sassumption} \label{epsmonotone}
The scheme $S$ is  $\epsilon(p,K)$ monotone.
\end{Sassumption}
\begin{Sassumption} \label{regular}
(Regularity of $S$ scheme) For every $h > 0$ and $\phi \in C_b(\R^N)$, $x \longrightarrow   S(h,x,\phi(x),[\phi]_x)$ is bounded and continuous in $\R^N$ and the function $r \longrightarrow S(h,x,r,[\phi]_x)$ is uniformly continuous for bounded $r$, uniformly in $x \in \R^N$.
\end{Sassumption}
\begin{Sassumption} \label{consistency}
(Consistency) There exists  a set of strictly positive integers  $(k_i)_{i =1, m}$, and a constant $K$ such that for every $h \ge 0$, $x \in \R^N$ and a smooth function $\phi$
\begin{eqnarray}
| F(x, \phi(x), \mathcal{D} \phi(x), \mathcal{D}^2 \phi) - S(h,x, \phi(x), [\phi]_x)| \le K  \sum_{i=1}^m h^{k_i} |\mathcal{D}^i \phi|_{0}. \nonumber
\end{eqnarray}
\end{Sassumption}
We add an assumption on the existence of a  solution of the discretized scheme with sufficient regularity which has to be checked for each scheme:
\begin{Sassumption}
\label{disregular}
We suppose there exists $C$ and $r$ independent of $h$ such that  for each $h$ we can find an  $\epsilon(C h^r)$ solution $u_h$ of scheme $S.$
\end{Sassumption}
\begin{remark}
Assumptions \ref{regular} and \ref{consistency} are classical. Assumption \ref{epsmonotone} is a relaxation of the constraint on the monotonicity of the scheme.\end{remark}
\begin{remark}
Related to assumption \ref{disregular}, \cite{bokanowski2015value} studied the case where some non monotone schemes were solved within
some margin error.
\end{remark}
An $\epsilon$ monotone scheme doesn't ensure  the existence of a discrete comparison result  but we give here a relaxation of this result:
\begin{lemma}
\label{discomp}
Assume \ref{epsmonotone}. Let $v \in C_b(\R^N)$ and $u \in C^{0,1}(\R^N)$. 
If $u$ is a subsolution of \reff{eqHJBDis} and $v$ is an  $\epsilon(C)$ supersolution  of \reff{eqHJBDis} then 
\begin{eqnarray}
u \le v +  \frac{ 1}{\bar \lambda}( K h^p |u|_1 + C). \nonumber
\end{eqnarray}
\end{lemma}
\proof
Mimicking  Lemma 2.3 in \cite{barles2002convergence}, we assume $m := \sup_{\R^N} (u-v) > \frac{1}{\bar \lambda} (K h^p |u|_1+  C)$.
Let $\{x_n\}_n$ be a sequence such that  $u(x_n) -v(x_n)  := \delta_n \longrightarrow  m $. For $n$ large enough, $\delta_n > \frac{1}{\bar \lambda} (K h^p |u|_1+  C)$.
Using the  subsolution definition, assumption \ref{epsmonotone} :
\begin{eqnarray}
0 & \ge & S(h, x_n, u(x_n), [u]_{x_n}) -S(h, x_n, v(x_n), [v]_{x_n})  -  C, \nonumber \\
0 &  \ge &  S(h, x_n, v(x_n) + \delta_n , [v + m ]_{x_n}) -S(h, x_n, v(x_n), [v]_{x_n})   -  K h^{ p} |u|_{1,0}  - C, \nonumber \\
0 &  \ge &   \bar \lambda \delta_n + w(m- \delta_n)  -  K h^{p} |u|_1  - C,   \nonumber 
\end{eqnarray}
where $w(t) \longrightarrow 0$ when $t^{+} \longrightarrow 0 $ by  assumption \ref{regular}. Letting $n \longrightarrow \infty$, we get
\begin{eqnarray}
m \le \frac{ 1}{\bar \lambda}\left( K h^p |u|_1+  C \right ),  \nonumber
\end{eqnarray}
which gives the contradiction.
\qedsymbol \\
\newline
We now give the lower and upper bound for the error given by the scheme. 
The lower bound will be given by the classical Krylov method of shaking coefficients \cite{krylov2000rate},
 while the upper bound will be given by the use of a switching system as in 
\cite{barles2005error,evans1979optimal} that gives a supersolution of the problem.
To use the theory developed in \cite{barles2005error}, we need to add one last assumption:
\begin{Aassumption}
\label{discommand}
For every $\delta > 0$, there are $Q \in \N$ and $\{\alpha_i\}_{i=1}^Q \subset \mathcal{A}$, such that for any $\alpha \in \mathcal{A}$
\begin{eqnarray}
\inf_{1 \le i  \le Q}  ( |\sigma^{\alpha} - \sigma^{\alpha_i}|_0 + | b^{\alpha} - b^{\alpha_i}|_0 + |c^{\alpha}- c^{\alpha_i}|_0 + |f^{\alpha} -f^{\alpha_i}|_0) < \delta.  \nonumber
\end{eqnarray}
\end{Aassumption}
We introduce the following switching system:
\begin{eqnarray}
\label{switching}
F_i^{\epsilon}(x,v^{\epsilon}, \mathcal{D} v^{\epsilon}_i ,\mathcal{D}^2 v^{\epsilon}_i) &= & 0  \mbox{ in } \R^N, i \in \mathcal{I} := \{ 1, .., Q \}, 
\end{eqnarray}
where $v^{\epsilon} = (v_1^{\epsilon}, ... , v^{\epsilon}_Q),$
\begin{eqnarray}
F_i^{\epsilon}(x,r,p,X)& =& max \left\{ \min_{|e| \le \epsilon} \mathcal{L}^{\alpha_i}(x + e,r_i,p,X); r_i - \mathcal{M}_i r \right \} , \nonumber
\end{eqnarray}
and $\mathcal{L}^{\alpha}$ given by \reff{Loperator}, 
\begin{eqnarray}
\mathcal{M}_i r & =&  \min_{j \ne i} \left \{ r_j + k \right \}.
\end{eqnarray}
We give  two lemmas proved in \cite{barles2005error}:
\begin{lemma}
\label{lemmeSwitch1}
Assume \ref{refAss1} and  \ref{refAss2}.
\begin{itemize}
\item There exists a unique solution $v^{\epsilon}$ of \reff{switching} satisfying $|v^{\epsilon}|_1 \le C$ where $C$ depends
only on $\lambda$, $K$ from \ref{refAss1}.
\item Assume in addition \ref{discommand}, then for any $\delta > 0$, there are $Q \in \N$ and $\{\alpha_i\}_{i=1}^Q \subset \mathcal{A}$ such that
the solution $v_{\epsilon}$ of \reff{switching} satisfies
\begin{eqnarray}
\max_i |u -v^{\epsilon}_i|_0 & \le & C (\epsilon + k^{\frac{1}{3}} + \delta),\nonumber
\end{eqnarray}
where $C$ depends on $\lambda$, $K$ from \ref{refAss1}.
\end{itemize}
\end{lemma}

A function $u$ can be regularized by 
\begin{eqnarray}
 \rho_{\epsilon} * u(x) := \int_{\R^N} u(x-e) \rho_{\epsilon}(e) de,   \nonumber \\
\end{eqnarray}
where $\rho_{\epsilon}$ is the mollifier sequence such that
\begin{eqnarray}
\rho_{\epsilon}(x) = \frac{1}{\epsilon^N} \rho(\frac{x}{\epsilon}) \mbox{ where } \rho \in C^{\infty}(\R^N),  \int_{\R^N} \rho = 1, \mbox{and supp}(\rho)  = \bar B(0,1). \nonumber
\end{eqnarray}

\begin{lemma}
\label{lemmeSwitch2}
Assume \ref{refAss1} and  \ref{refAss2} and define $v_{\epsilon,i} := \rho_{\epsilon} * v^{\epsilon}_i$ for $i \in \mathcal{I}$.
\begin{itemize}
\item There is a constant $C$ depending only on $\lambda$, $K$ from \ref{refAss1}, such that 
\begin{eqnarray}
|v_{\epsilon_,j} - v^{\epsilon}_i|_0 \le C( k + \epsilon) \mbox{ for } i,j \in \mathcal{I}. \nonumber
\end{eqnarray}
\item Assume in addition that $\epsilon  \le (4 \sup_i[v_i^{\epsilon}]_1) ^{-1}k$. For every $x \in \R^N$, if $j := argmin_{i \in \mathcal{I}} v_{\epsilon,i}(x)$,
then
\begin{eqnarray}
\mathcal{L}^{\alpha_j}(x, v_{\epsilon,j}(x), \mathcal{D} v_{\epsilon,j}(x), \mathcal{D}^2 v_{\epsilon,j}(x)) \ge 0.
\end{eqnarray}
\end{itemize}
\end{lemma}
We can now give the main result of the paper using the two previous lemmas for an  upper bound of the $\epsilon$ solution.
\begin{theorem}
\label{mainTheo}
Assume \ref{refAss1},  \ref{refAss2}, \ref{epsmonotone}, \ref{regular}, \ref{consistency}, and \ref{disregular}.
We have the following bounds for $u_h$ a sequence of $\epsilon(\tilde Ch^r)$ solutions:
\begin{eqnarray}
u-u_h & \le &  \hat C(h^{\min(p,r)} +h^{\underset{i=1,m}{\min}\frac{k_i}{i}}), \nonumber
\end{eqnarray}
where $\hat C$ depends on $\tilde C$.\\
Besides assume  \ref{discommand} then there exists $\hat C$ depending on $\tilde C$ such that :
\begin{eqnarray}
u_h -u & \le & \hat C ( h^{\min(p,r)} + h^{\underset{i=1,m}{\min} \frac{k_i}{3i-2}}). \nonumber 
\end{eqnarray}
\end{theorem}
\proof
For the lower bound, we follow the Krylov demonstration  \cite{krylov2000rate} as  done in \cite{barles2002convergence}.
First we introduce  the solution $u^{\epsilon}$
\begin{eqnarray}
\max_{|e| \le \epsilon} [ F(x+e, u^{\epsilon}, \mathcal{D} u^{\epsilon},  \mathcal{D}^2 u^{\epsilon}) ] = 0  \mbox{ in } \R^N.  \nonumber
\end{eqnarray} 
The existence a solution $u^{\epsilon}$ in $C^{0,1}(\R^N)$  satisfying $|u^{\epsilon}| \le C$ and $|u^{\epsilon} -u |_0 \le C \epsilon$ is
given by Lemma 2.6  in \cite{barles2002convergence}.
Noting that for each $e \le \epsilon$, $u^{\epsilon}(.-e)$ satisfies for each function $\phi= \psi(.-e) \in C^2(\R^N)$ and each $y$ where $u^{\epsilon}(y-e) - \phi(y)$ is maximal
\begin{flalign*}
F( y, u^{\epsilon}(y-e), \mathcal{D} \phi(y), \mathcal{D}^2 \phi(y))   = & F( y, u^{\epsilon}(y-e), \mathcal{D} \psi(y-e), \mathcal{D}^2 \psi(y-e)),  \nonumber \\
 \le & \sup_{|e| \le \epsilon} F( y, u^{\epsilon}(y-e), \mathcal{D} \psi(y-e), \mathcal{D}^2 \psi(y-e)), \nonumber\\
 \le & 0,  \nonumber
\end{flalign*}
so  $u^{\epsilon}(.-e)$ is a subsolution of \reff{eqHJB}.\\
Then $u^{\epsilon}(.-e)$ is regularized by
\begin{eqnarray}
u_{\epsilon}(x) := \int_{\R^N} u^{\epsilon}(x-e) \rho_{\epsilon}(e) de. \nonumber
\end{eqnarray}
The regularized function $u_{\epsilon}$ is a subsolution of problem \reff{eqHJB} as given by Lemma 2.7 in \cite{barles2002convergence}.
First use the relation  for $m>0$, 
$$F(x, t+m, p ,X) \ge F(x,t ,p ,X) + \lambda m,$$
 and the consistency property \ref{consistency} to get
\begin{flalign*}
F(y, u_{\epsilon}(y), \mathcal{D} u_{\epsilon}(y), \mathcal{D}^2  u_{\epsilon}(y) )  \ge & F(y,u_{\epsilon}(y)- \frac{K}{\lambda} \sum_{i=1}^m  h^{k_i} |\mathcal{D}^i u_{\epsilon}|_{0}, \mathcal{D} u_{\epsilon}(y), \mathcal{D}^2  u_{\epsilon}(y))  + \\
&  K \sum_{i=1}^m  h^{k_i} |\mathcal{D}^i u_{\epsilon}|_{0},   \\
 \ge & S(h, y, u_{\epsilon}(y)- \frac{K}{\lambda} \sum_{i=1}^m  h^{k_i} |\mathcal{D}^i u_{\epsilon}|_{0}, [u_{\epsilon}- \frac{K}{\lambda} \sum_{i=1}^m  h^{k_i} |\mathcal{D}^i u_{\epsilon}|_{0}]_y). 
\end{flalign*}
Then $u_{\epsilon}-  \frac{K}{\lambda} \sum_{i=1}^m  |\mathcal{D}^i u_{\epsilon}|_{0}$ is a subsolution of equation \reff{eqHJBDis}.
Using lemma \ref{discomp}, and  assumption \ref{disregular}, we get that there exists $C$  such that
\begin{eqnarray}
u_{\epsilon} - u_h  & \le  & C (  |u_{\epsilon}|_1  h^{p} +   \sum_{i=1}^m  |\mathcal{D}^i u_{\epsilon}|_{0}  h^{k_i} + h^r), \nonumber \\
  & \le & C ( h^{\min(p,r)} + \sum_{i=1}^m \epsilon^{1-i} h^{ k_i}).  \nonumber
\end{eqnarray}
In the last line we have used that because  $u^{\epsilon}$ is bounded uniformly in $C^{0,1}(\R^N)$, $u_{\epsilon}$ is regular  and $|\mathcal{D}^n u_{\epsilon}|_0 \le C \epsilon^{1-n}$ for $n \ge 1$.\\
At last using the mollifier properties, the uniform boundedness of $u^{\epsilon}$ in $C^{0,1}$ gives us that $|u_{\epsilon}- u^{\epsilon}| \le C \epsilon$.
Besides $|u^{\epsilon} -u |_0 \le C \epsilon$ (Lemma 2.6 in \cite{barles2002convergence})  so there exists $\hat C$ depending on $\tilde C$ such that
\begin{eqnarray}
u - u_h  & \le & |u-u^{\epsilon}|_0 +  |u^{\epsilon}-u_{\epsilon}|_0 + u_{\epsilon} - u_h, \nonumber \\
         & \le &  C (h^{\min(p,r)}+  \sum_{i=1}^m \epsilon^{1-i} h^{ k_i}+ \epsilon). \nonumber 
\end{eqnarray}
Choosing $\epsilon = h^{\underset{i=1,m}{\min} \frac{k_i}{i}}$ we get the lower bound in the theorem.\\
For the upper bound we follow the Barles Jakobsen demonstration that we shorten except for points different from the initial proof.
We fix a $\delta > 0$, and pick up the corresponding $\{\alpha_i\}_{i\in \mathcal{I}}$ according to \ref{discommand}.
The corresponding solution $v^{\epsilon}$ of \reff{switching} exists according to lemma \ref{lemmeSwitch1} and is regularized as in lemma  \ref{lemmeSwitch2}.
We note 
\begin{eqnarray}
m := \sup_{y \in \R^N} \{ u_h(y) - w(y) \}, \nonumber
\end{eqnarray}
where $ w := \min_{i \in \mathcal{I}} v_{\epsilon,i}.$ 
We approximate $m$ by 
\begin{eqnarray}
\label{mkappa}
m_{\kappa} := \sup_{y \in \R^N} \{ u_h(y) - w(y) - \kappa \phi(y) \},
\end{eqnarray}
where $\phi(y)= (1+|y|^2)^{\frac{1}{2}}$.
Because of the boundedness and continuity of $u_h$ and $w$, the maximum is attained at a point x.
Because of the definition of $x$,
\begin{eqnarray}
m_{\kappa} := \sup_{y \in \R^N} \{ u_h(y) - v_{\epsilon,i}(y) - \kappa \phi(y) \}, \nonumber
\end{eqnarray}
where $i = argmin_{i \in \mathcal{I}} v_{\epsilon,i}(x)$.\\
Taking $\epsilon  = (4 \sup_i[v_i^{\epsilon}]_1) ^{-1}k$, from lemma \ref{lemmeSwitch2}, the definition of $\phi$ and \ref{refAss1} we get
\begin{eqnarray} 
\sup_{\alpha \in \mathcal{A}} \mathcal{L}^{\alpha}(x, (v_{\epsilon_,i}+ \kappa \phi)(x),\mathcal{D} (v_{\epsilon_,i}+ \kappa \phi)(x), \mathcal{D}^2 (v_{\epsilon_,i}+ \kappa \phi)(x)) \ge -C \kappa. \nonumber
\end{eqnarray}
Then using \ref{consistency}:
\begin{eqnarray}
-C  \kappa & \le &  S(h, (v_{\epsilon_,i}+ \kappa \phi)(x), [v_{\epsilon_,i}+ \kappa \phi]_x) + K \sum_{j=1}^m  h^{k_j} | \mathcal{D}^j (v_{\epsilon_,i}+ \kappa \phi)|_{0}. \nonumber
\end{eqnarray}
Using the properties of the mollified $v_{\epsilon,i}$, and the definition of $\phi$ :
\begin{eqnarray}
     -K  \sum_{j=1}^m  h^{k_j} \epsilon^{1-j}+  \mathcal{O}(\kappa)  &\le &  S(h, (v_{\epsilon_,i}+ \kappa \phi)(x), [v_{\epsilon_,i}+ \kappa \phi]_x). \label{mKappa1}
\end{eqnarray}
Then we use the $\epsilon$ monotony property \reff{monoton2}, the definition of $m_{\kappa}$, the fact that $v_{\epsilon_,i}$ is bounded uniformly by the properties of mollifiers and lemma \ref{lemmeSwitch1}
  to get
\begin{eqnarray}
S(h, (v_{\epsilon_,i}+ \kappa \phi)(x), [v_{\epsilon_,i}+ \kappa \phi]_x) & \le & S(h, u_h(x)- m_{\kappa}, [ u_h- m_{\kappa}]_x)  + K h^p |v_{\epsilon_,i}+ \kappa \phi|_1, \nonumber\\
 & \le &  - \lambda m_{\kappa} + S(h, u_h(x), [ u_h]_x)  + C h^p (1+ \kappa), \nonumber \\
& \le &   \tilde C (1+ \kappa) h^{ \min(r,p)} - \lambda m_{\kappa}.  \label{mKappa2}
\end{eqnarray}
Using \reff{mKappa1} and \reff{mKappa2} we get
\begin{eqnarray}
\lambda m_{\kappa} & \le &    \tilde C (1+ \kappa) h^{ \min(r,p)}  + K \sum_{j=1}^m  h^{k_j} \epsilon^{1-j} + \mathcal{O}(\kappa). \nonumber
\end{eqnarray}
An estimate of $m$ is obtained by letting $\kappa$ goes to $0$.
Then for any $y \in \R^N$,
\begin{eqnarray}
u_h(y) -u(y) & \le & u_h(y) - v_{\epsilon,i}(y) + v_{\epsilon,i}(y) - u(y), \nonumber \\
             & \le & m +  v_{\epsilon,i}(y) - u(y). \nonumber
\end{eqnarray}
Using the lemma \ref{lemmeSwitch1} and \ref{lemmeSwitch2} we get
\begin{eqnarray}
u_h(y) -u(y) & \le &  \hat C(  h^{ \min(r,p)}   + \sum_{i=1}^m  h^{k_i} \epsilon^{1-i} + \epsilon + k + k^{\frac{1}{3}} + \delta),  \nonumber
\end{eqnarray}
with $\hat C$  depending on  $\tilde C$ and uniform in $y$.
The conclusion is obtained by taking $ \epsilon = h^{\underset{i=1,m}{\min} \frac{3k_i}{3i-2}}$ , remembering that $k=\epsilon \sup_i[v_i^{\epsilon}]_1$  and letting $\delta$ going to 0.
\qedsymbol\\
\section{Some numerical schemes}
In this  section  we take the notations and we will follow arguments similar to \cite{barles2002convergence,barles2005error}. We give the notations used for the discretization and interpolation used by the scheme.
A thorough study of interpolation method for time dependent HJB equation can be found in \cite{warin2016some}. All the schemes defined in \cite{warin2016some} can be used :
it includes truncated Lagrangian interpolators, classical cubic spline truncated interpolators, the monotone cubic spline first defined in \cite{debrabant2013semi}.
In the sequel we focus on the  truncated  GLL interpolators which are the most effective according to \cite{warin2016some}.

A spatial discretization $\Delta x$  of the problem being given,  in the sequel a mesh  $\hat X_{\bar i}$ corresponds to the hyper-cube
  $[i_1 \Delta x , (i_1+1) \Delta x] \times ... \times [i_N \Delta x, (i_N+1) \Delta x]$  with  $\bar i = (i_1 , ..., i_N) \in \mathbf{ Z}^N$ .
For GLL quadrature grid $(\xi_{i})_{i=1,...,M+1} \in [-1,1]$, with $\xi_1=-1$, $\xi_{M+1}=1$, and for a mesh $\bar i$, the point $y_{\bar i, \tilde j}$ with $\tilde j = (j_1, ..., j_N) \in [1,M+1]^N$  will have the coordinate
$(\Delta x (i_1 + 0.5 (1+\xi_{j_1})) ,...,\Delta x (i_N + 0.5 (1+\xi_{j_N}))$.
We denote  $ X_{\Delta x, M} := (y_{\bar i, \tilde j})_{\bar i, \tilde j}$ the set of all the grids points on the whole domain.\\
We notice that  each mesh $\hat X_{\bar i}$ has  a constant volume  $\Delta x^N$, so we have the following relation for all $x \in \R^N$:
\begin{eqnarray}
\min_{\bar i, \tilde j} | x- y_{\bar i, \tilde j}| \le  C \Delta x.
\label{MaxDist}
\end{eqnarray}
We introduce $I_{\Delta x,M}$ the  Lagrange interpolator associated to the   GLL quadrature. 
We recall that in one dimension, the GLL Lagrange interpolator $I_M$ on $[-1,1]$
is given by (see for example  \cite{quarteroni2010numerical}):
\begin{eqnarray*}
I_{M}(f) & =&  \sum_{k=0}^{M} \tilde f_k L_k(x) , \nonumber  \\
\tilde f_k &  =  &   \frac{1}{\gamma_k} \sum_{i=0}^{M} \rho_i f(\eta_i) L_k( \eta_i) , \nonumber \\
\gamma_k & = & \sum_{i=0}^{M} L_k(\eta_i)^2 \rho_i,
\end{eqnarray*}
where the functions  $L_N$  satisfy the recurrence
\begin{eqnarray*}
(N+1) L_{N+1}(x) & = &  (2N+1) x L_N(x) -N L_{N-1}(x), \\
L_0 & = & 1 , \quad   L_1 = x , 
\end{eqnarray*}
$\eta_1 = -1, \eta_{M+1} =1$, the $\eta_i$ $(i=2,...,M)$  are the  zeros of $L^{'}_{M}$  and the eigenvalues of the matrix $P$
\begin{eqnarray*}
P &=& \left ( \begin{array}{lllll}
              0  & \gamma_1 & ...& 0 & 0 \\
              \gamma_1 & 0 & ...   & 0 & 0 \\
              ...      & ... & ... & ... & ... \\
              0  & 0 & ... & 0 & \gamma_{M-2} \\
              0  & 0 & ... & \gamma_{M-2} & 0 
              \end{array}
 \right), \nonumber  \\
\gamma_n & = & \frac{1}{2} \sqrt{\frac{n(n+2)}{(n+\frac{1}{2})(n+\frac{3}{2})}}  , 1 \le n \le M-2 , \nonumber
\end{eqnarray*}
and the weights satisfies 
\begin{eqnarray*}
\rho_{i} &  = & \frac{2.}{(M+1)M L_{M}^2(\eta_i)} ,  1 \le i \le M+1.
\end{eqnarray*}
The interpolator $I_{\Delta x , M}$ on a mesh $\hat X_{\bar i}$  is obtained by first rescaling $I_M$   and by  tensorization.

On a mesh $\hat X_{\bar i}$, we note  $\underline v_{\bar i} = \displaystyle  \min_{\tilde j} v(y_{\bar i, \tilde j})$, $\bar v_{\bar i} = \displaystyle  \max_{\tilde j} v(y_{\bar i, \tilde j})$. 
We introduce the following truncated operator:
\begin{eqnarray}
\hat I_{\Delta x, M }(v) & =&  \underline v_{\bar i}   \vee I_{\Delta x, M }(v) \wedge  \bar v_{\bar i},  \nonumber 
\end{eqnarray}
where $\wedge$ denotes the minimum and $ \vee$ the maximum.\\
We first give some properties associated to the truncated interpolation operator.
\begin{lemma}
\label{weights}
The interpolator $\hat I_{\Delta x, M }$ satisfies 
\begin{eqnarray}
(\hat I_{\Delta x, M} f)(x) & = & \sum_{\tilde j} (w_{\bar i, \tilde j}(f))(x) f(y_{\bar i, \tilde j})  \nonumber \\
 \sum_{\tilde j} (w^h_{\bar i, \tilde j}(f))(x)  & =& 1, \nonumber
\end{eqnarray} 
and the positive weights $w^h_{\bar i, \tilde j}(f)$ are functions of $f$. 
\end{lemma}
\proof
Because of the truncation for each point  $x$ of a mesh $\hat X_{\bar i}$, we have 
\begin{eqnarray}
(\hat  I_{\Delta x, M} (f)(x)  =  \underline w^h_{\bar i}(f)(x) \underline f_{\bar i} + \bar w^h_{\bar i}(f)(x) \bar  f_{\bar i}. \nonumber
\end{eqnarray}
where the  weights are defined as follows:\\
If $ \underline f_{\bar i}  \le  \hat  I_{\Delta x, M} (f)(x) \le  \bar  f_{\bar i}$ then
\begin{eqnarray}
 \underline w^h_{\bar i}(f)(x) & = & \frac{\hat  I_{\Delta x, M}(f)(x) - \bar  f_{\bar i}}{ \underline f_{\bar i} - \bar  f_{\bar i}},  \nonumber \\
\bar w^h_{\bar i}(f)(x) & = & 1 -\underline w^h_{\bar i}(f)(x),  \nonumber\\
\end{eqnarray}
If $ \underline f_{\bar i} > \hat  I_{\Delta x, M} (f)(x)$,
\begin{eqnarray}
 \bar w^h_{\bar i}(f)(x) &=& 0, \nonumber \\
\underline w^h_{\bar i}(f)(x)& =&  1.  \nonumber 
\end{eqnarray}
Otherwise
\begin{eqnarray}
 \underline w^h_{\bar i}(f)(x) &=& 0, \nonumber  \\
\bar w^h_{\bar i}(f)(x)& =& 1.  \nonumber
\end{eqnarray}
The weights associated to non extremal points are taken equal to $0$.
\qedsymbol\\
Then we  add a result for the interpolation error :
\begin{lemma}
 \label{interpolerror}
\begin{itemize}
\item 
\label{interpolerrorLip}
For each K-Lipschitz bounded function $f$:
\begin{eqnarray}
| \hat I_{\Delta x, M}(f) -f |_0 & \le &  K (\Delta x ).  \nonumber
\end{eqnarray}
\item 
\label{interpolerrorReg}
Suppose $M \ge 2$ , for each  $f$  $\R$ value function defined on $\R^N$ and twice differentiable , there exists  $C$ such that:
\begin{eqnarray}
| \hat I_{\Delta x, M}(f) -f |_0 & \le &  C  \Delta x^2 |D^2f|_{0}.  \nonumber
\end{eqnarray}
\item
\label{interpolLinear}
For $m \in R^N$, 
\begin{eqnarray}
 \hat I_{\Delta x, M}(f+m)(x) =  \hat I_{\Delta x, M}(f)+m. \nonumber
\end{eqnarray}
\end{itemize}
\end{lemma}
\proof
First because of the truncation, continuity of $\hat I_{\Delta x, M}(f)$, for each $x \in \R^N$, $x \in \hat X_{\bar i}$, 
 there exists $\tilde x \in  \hat X_{\bar i}$ such that $\hat I_{\Delta x, M}(f)(x)=f(\tilde x)$. We then use the Lipschitz property of
$f$ and  \reff{MaxDist} to get the result.\\
When no truncation is achieved, we have a rate of convergence in $O(\Delta x^{M+1})$. 
When the truncation is achieved, for example $\hat I_{\Delta x, M}(f)(x)= \bar f_{\bar i} $, then
\begin{eqnarray*}
 I_{\Delta x,1} \le \hat I_{\Delta x, M}(f)(x)  \le I_{\Delta x, M}(f)(x),
\end{eqnarray*}
where $I_{\Delta x,1}$ correspond to the linear interpolator and then
\begin{eqnarray*}
| \hat I_{\Delta x, M}(f)(x) -f(x)| \le \max ( |I_{\Delta x, M}(f)(x) -f(x)|, |I_{\Delta x, 1}(f)(x) -f(x)|),
\end{eqnarray*}
so the rate of convergence remains at least equal to 2.\\
The third point is easily check by noticing that $I_{\Delta x, M}$ is a Lagrange interpolator so linear, that  $I_{\Delta x, M}(m) = m$ and that
the truncation operator $tr$ satisfies $tr(f+m) = tr(f) +m$.
\qedsymbol\\
\begin{remark}
Some effective interpolation methods such as ENO, WENO can be used  for interpolation \cite{jiang2000weighted,osher1991high,shu2007high} while solving Hamilton Jacobi equations but they are not proved convergent. 
We will show that the previously defined interpolator
ensures convergence for Semi Lagrangian schemes and for some Finite Difference schemes  but at a rate not better than linear interpolator.
The interest of such interpolators will by checked numerically on some examples  in the last section. Besides, they are easy to implement independently of the dimension of the problem.
\end{remark}
\subsection{A Camilli Falcone style scheme}
The first scheme we study is a modification of the scheme of Camilli and  Falcone  \cite{camilli1995approximation} where the linear interpolator $I_{\Delta x,1}$ is replaced by a potentially high order
interpolator $\hat I_{\Delta x,M}$ with $M>1$.
We begin by defining the monotone operator \cite{barles2002convergence,camilli1995approximation} $\hat S$ which is the Lagrangian scheme without interpolation.
First for any bounded continuous function $\phi$, $x,y,z \in \R^N$, we set 
\begin{eqnarray}
\label{scheme0}
\hat S(h,y,t,\phi_x) =  \sup_{\alpha \in \mathcal{A}} \left \{ -\frac{1}{h}(G(h,\alpha,y,\phi_x)-t) + c^{\alpha}(y) t - f^{\alpha}(y) \right \} , \nonumber
\end{eqnarray} 
\begin{eqnarray}
G(h,\alpha,y,\phi_x) & = & \frac{1-h c^{\alpha}(y)}{2P} \sum_{i=1}^P \left( \phi(x+h b^{\alpha}(y) + \sqrt{hP} \sigma^{\alpha}_{i}(y)) + \phi(x+h b^{\alpha}(y) -
 \sqrt{hP} \sigma^{\alpha}_i(y)) \right ), \nonumber
\end{eqnarray}
where $\sigma^{\alpha}_{i}$ is the i-th column of $\sigma$.\\
The semi discretized scheme $S$ is defined as follows  for $y$ a quadrature point.
\begin{eqnarray}
\label{scheme1}
S(h,y,t, [\phi]_x) =  \hat S(h, y, t,  (\hat I_{\Delta x ,M} \phi)_x). \nonumber
\end{eqnarray}
So the discretized problem leads to find $U$ function on $X_{\Delta x ,M}$ such that
\begin{eqnarray}
\label{discreteHJB}
| \hat S(h, y, U(y),  (\hat I_{\Delta x ,M} U)_y)|   & \le & \epsilon(h, \Delta x) , \mbox{ for } y \in X_{\Delta x,M}.
\end{eqnarray}
We first recall some results on the solution associated to the semi discretized scheme that can be found in \cite{camilli1995approximation,barles2002convergence}
\begin{proposition}
Assume that \ref{refAss1}, \ref{refAss2} hold. Then there exists a unique bounded function $v_h$  uniformly  in $C^{0,1}(\R^N)$ satisfying
\begin{eqnarray}
\label{scheme00}
 \hat S(h,x,v_h(x), (v_h)_x) &= & 0 \quad \mbox{ for } x \in \R^N.
\end{eqnarray}
\end{proposition} 
We next prove that the scheme  $S$ satisfies the first assumptions of the article :
\begin{proposition}
Assume \ref{refAss1} hold and that $\Delta x= h^{q}$. 
Then the scheme \reff{scheme1} satisfies assumptions \ref{epsmonotone}, \ref{regular}, \ref{consistency} with
$k_2 = \min(2q-1,1)$, $k_4= 1$ , $p=q-1$.
\end{proposition}
\proof
First assumption \ref{regular} follows easily from \ref{refAss1}. 
\ref{consistency} follows easily using  lemma \ref{interpolerror}: for $v$ regular  the consistency error is bounded  by
\begin{eqnarray}
\epsilon(h,\Delta x) & = C ( h |\mathcal{D}^4 v |_0 + h |\mathcal{D}^2 v |_0 + \frac{\Delta x^2}{h} |\mathcal{D}^2 v |_0).
\end{eqnarray}
Using  \ref{refAss1}, and  lemma \ref{interpolerror}, for $m > 0$ :
\begin{eqnarray}
G(h,\alpha,y,(\hat I_{\Delta x ,M}(\phi+m))_x) & =  &  G(h,\alpha,y,(\hat I_{\Delta x ,M}(\phi))_x) + (1 - h c^{\alpha}(x)) m , \nonumber
\end{eqnarray}
so
\begin{eqnarray}
S(h,y,t+m, [ I_{\Delta x ,M}(\phi+m)]_x) & \ge    &  S(h,y,t, [ I_{\Delta x ,M}(\phi)]_x)  + 2 \lambda m, \nonumber
\end{eqnarray}
and property \reff{monoton2} is checked.\\
Suppose $v\in C_b(\R^N)$ , $w \in C^{0,1}(\R^N)$,  $v \ge w$.
If $x \in \hat X_{\bar i}$ is such that $\hat I_{\Delta x ,M}(v)(x) \le \hat I_{\Delta x ,M}(w)(x)$ let's introduce  $v(x_{\bar i, \tilde l}) = \min_{\tilde j} v(x_{\bar i, \tilde j})$.
It satisfies  $v(x_{\bar i, \tilde l}) \le \hat I_{\Delta x ,M}(v)(x)$ so using lemma \ref{interpolerror}
\begin{eqnarray}
\hat I_{\Delta x ,M}(v)(x) -\hat I_{\Delta x ,M}(w)(x) \ge v(x_{\bar i, \tilde l}) - \hat I_{\Delta x ,M}(w)(x),  \nonumber \\
                                       \ge  v(x_{\bar i, \tilde l})  - w(x_{\bar i, \tilde l}) - |w|_1 \Delta x \sqrt{N}, \nonumber
                                       \ge  - |w|_1 \Delta x  \sqrt{N}.
\end{eqnarray}
So for all $x \in \R^N$
\begin{eqnarray}
\hat I_{\Delta x ,M}(v)(x) & \ge & \hat I_{\Delta x ,M}(w)(x)  - |w|_1 \Delta x  \sqrt{N},
\label{estimEsp1}
\end{eqnarray}
and 
\begin{eqnarray}
G(h,\alpha,y,(\hat I_{\Delta x ,M}(v))_x) & \ge   &  G(h,\alpha,y,(\hat I_{\Delta x ,M}(w))_x) - (1 + h |c^{\alpha}(y)|_0) |w|_1 \Delta x  \sqrt{N}. \nonumber
\end{eqnarray}
So
\begin{eqnarray}
S(h,y,t,[ \hat I_{\Delta x ,M}(v)]_x) \le   S(h,y,t,[ \hat I_{\Delta x ,M}(w)]_x)  +  (1 + h |c^{\alpha}(y)|_0) |w|_1 \frac{ \sqrt{N}\Delta x}{h}.  \nonumber
\end{eqnarray}
Similarly if $v\in C^{0,1}(\R^N)$ , $w \in C_b(\R^N)$,  $v \ge w$, noting $w(x_{\bar i, \tilde l}) = \max_{\bar j} v(x_{\bar i, \tilde j})$,
and using lemma \ref{interpolerror}
\begin{eqnarray}
\hat I_{\Delta x ,M}(v)(x) -\hat I_{\Delta x ,M}(w)(x) \ge \hat I_{\Delta x ,M}(v)(x) -w(x_{\bar i, \tilde l}), \nonumber \\
                                       \ge  v(x_{\bar i, \tilde l})  - w(x_{\bar i, \tilde l}) - |v|_1 \Delta x  \sqrt{N} \nonumber
                                       \ge  - |v|_1 \Delta x  \sqrt{N},
\end{eqnarray} 
so
\begin{eqnarray}
\hat I_{\Delta x ,M}(v)(x) & \ge     
                   & \hat I_{\Delta x ,M}(w)(x)  - |v|_1 \Delta x  \sqrt{N},
\label{BestimEsp2}
\end{eqnarray}
and
\begin{eqnarray}
S(h,y,t,[ \hat I_{\Delta x ,M}(v)]_x) \le   S(h,y,t,[ \hat I_{\Delta x ,M}(w)]_x)  +   (1 + h |c^{\alpha}(y)|_0)  |v|_1 \frac{ \sqrt{N} \Delta x}{h}, \nonumber
\end{eqnarray}
so that  the \reff{monoton1} and \reff{monoton1bis} properties are checked.
\qedsymbol\\

We need to prove that we can construct an approximate solution of the discretized problem.
We introduce the operator $T$ defined for  $U$ a function on $X_{\Delta x, M}$ :
\begin{flalign*}
(T_{h,\Delta x} U)(x) =&  \inf_{\alpha \in \mathcal{A}} \left \{ (1- h c^{\alpha}(x)) (\Pi_{\Delta x,h}(U))(x)  + h f^{\alpha}(x) \right \}  \mbox{ for } x \in X_{\Delta x, M}, \nonumber
\end{flalign*}
where the operator $\Pi_{\Delta x,h}$ is
\begin{eqnarray}
(\Pi_{\Delta x,h} U)(x) &= & \frac{1}{2P} \sum_{i=1}^{2P} \left( (\hat I_{\Delta x ,M} U)(x+ h b^{\alpha}(x) + \sqrt{Ph} \sigma^{\alpha}_i(x)) +(\hat I_{\Delta x ,M} U)(x+ h b^{\alpha}(x) - \sqrt{Ph} 
\sigma^{\alpha}_i(x)) \right). \nonumber
\end{eqnarray}
We then recursively define $T_{h,\Delta x}^s$ for $s \in \N$  and $s \ge 2$ by
\begin{eqnarray}
(T_{h,\Delta x}^s U)(x) = (T_{h,\Delta x} (T_{h,\Delta x}^{s-1} U))(x). 
\label{defineTs}
\end{eqnarray}
\begin{proposition}
\label{epsSolution1}
Assume \ref{refAss1}, \ref{refAss2} hold. Suppose that   $\Delta x= h^{q}$ with $q>2$.
There exists $s \in \N$ depending on $h$ and  $C$ independent of $h$ such that  $u_h =  \hat I_{\Delta x,M}(T_{h,\Delta x}^s 0)$ is an $\epsilon(C h^{q-2})$ solution of scheme $S.$
\end{proposition}
\proof
Note $v_h$ the unique solution given of scheme \reff{scheme00}.\\
For $x \in X_{\Delta x, M}$, $U$ a function on $X_{\Delta x, M}$, using  $| \inf A - \inf B| \le \sup | A - B|$
\begin{flalign}  \label{diffth}
|(T_{h,\Delta x} U)(x) -v_h(x)|  \le &  (1 -\lambda h)  \sum_{i=1}^{2P} \left(  \sup_{\alpha \in \mathcal{A}}  |(\hat I_{\Delta x ,M} U)(x+ h b^{\alpha}(x) +\sqrt{hP} \sigma^{\alpha}_i(x)) -  \right. \nonumber \\
& 
v_h(x+ h b^{\alpha}(x) +\sqrt{hP} \sigma^{\alpha}_i(x))|  \nonumber\\
  &   \left .  + |(\hat I_{\Delta x ,M} U)(x+ h b^{\alpha}(x) -\sqrt{hP} \sigma^{\alpha}_i(x)) 
  - v_h(x+ h b^{\alpha}(x) +\sqrt{hP} \sigma^{\alpha}_i(x))|  \right), 
\end{flalign}
then notice that for $x \in \hat X_{\bar i}$
\begin{eqnarray}
|\hat I_{\Delta x ,M} (U)(x) -v_h(x)| & = & | \sum_{\tilde j} w_{\bar i, \tilde j}(U)(x) (U(x_{\bar i, \tilde j}) - v_h(x)) |, \nonumber \\
& \le &  \sum_{j,l} w_{\bar i, \tilde j}(U)(x)  | U(x_{\bar i, \tilde j}) - v_h(x)|,  \nonumber \\
& \le & \sup_{\tilde j} | U(x_{\bar i, \tilde j}) - v_h(x)|,  \nonumber \\
& \le & \sup_{\tilde j}  | U(x_{\bar i, \tilde j}) - v_h(x_{\bar i, \tilde j})| + \Delta x  \sqrt{N} |v_h|_1, \nonumber \\
& \le & |U- v_h|_0 + C  \Delta x, \label{interpest}
\end{eqnarray}
where we have use the uniform boundedness of $v_h$ in $C^{0,1}(\R^N)$ .
Gathering \reff{diffth} and \reff{interpest} we get
\begin{eqnarray}
|(T_{h,\Delta x} U)(x) -v_h(x)| & \le & (1- \lambda h) (|U- v_h|_0 + C \Delta x  ). \nonumber
\end{eqnarray}
Iterating we find that
\begin{eqnarray}
|(T_{h,\Delta x}^k 0)(x) -v_h(x)| & \le & (1- \lambda h)^k |v_h|_0 + C \frac{\Delta x}{h}. 
\end{eqnarray}
Taking $k= \min( i \in \N \mbox{ such that } i \ge \frac{(q-1) log(h)}{log(1- \lambda h)})$ , using the fact that $\Delta x = h^{q}$, we get that
\begin{eqnarray}
|(T_{h,\Delta x}^k 0)(x) -v_h(x)|   \le C h^{q-1}. \label{iterationT}
\end{eqnarray}
Let's prove that   $u_h = \hat I_{\Delta x,N} (T_{h,\Delta x}^k 0)$ is an $\epsilon(h^{q-2})$ solution of the scheme $S$ .\\
As in  \reff{interpest}
\begin{eqnarray}
|\hat I_{\Delta x ,M} (T_{h,\Delta x}^k 0)(x) -v_h(x)| & \le & |T_{h,\Delta x}^k 0- v_h|_0 + C  \Delta x, \nonumber \\
& \le &  C h^{q-1}, \nonumber 
\end{eqnarray}
where we have used  \reff{iterationT}, so
\begin{eqnarray}
|u_h(x) -v_h(x)| \le   C h^{q-1} \label{directBound}.
\end{eqnarray}
Then using \ref{refAss1}, the fact that $u_h= \hat I_{h,\Delta x} u_h$, and \reff{directBound}
\begin{eqnarray}
| S(h,x,u_h(x), [u_h]_x)|  & =  & | \hat S(h,x,u_h(x), (u_h)_x) - \hat S(h,x,v_h(x), (v_h)_x)|,  \nonumber \\
                          & \le  &   \sup_{\alpha \in \mathcal{A}}  |\frac{1}{h} [ \left( G(h,\alpha,x,(u_h)_x)-u_h(x) \right) -  \left( G(h,\alpha,x,(v_h)_x)-v_h(x) \right) ]|,  \nonumber \\
                          & \le & c h^{q-2}. \nonumber
\end{eqnarray}
\qedsymbol\\

\begin{proposition}
\label{SemiLagProp1}
Assume \ref{refAss1}, \ref{refAss2}, \ref{discommand}  hold. Suppose $u_h$ has being constructed as in proposition \ref{epsSolution1}, the previous developed framework gives us that we can find $q$ such  that 
\begin{eqnarray}
 |u -u_h| < C h^{\frac{1}{10}}. \nonumber 
\end{eqnarray}
\end{proposition}
\begin{proposition}
\label{SemiLagProp2}
Assume \ref{refAss1}, \ref{refAss2}  hold. Taking $q \ge \frac{5}{4}$,   $u_h$  being constructed as in proposition \ref{epsSolution1}, we have
\begin{eqnarray}
 |u -u_h| < C h^{\frac{1}{4}}. \nonumber 
\end{eqnarray}
\end{proposition}
\proof
Because  we have a bound on $|v_h-u|_0$  in $O(h^{\frac{1}{4}})$ (see \cite{barles2002convergence}), a direct use of \reff{directBound} shows 
\begin{eqnarray}
| u- u_h| & \le & |u- v_h|_0 + |v_h- u_h|_0, \nonumber \\
          & \le & c (h^{\frac{1}{4}} + h^{q-1}),      \nonumber \\
\end{eqnarray}
giving the result.
\qedsymbol\\

\begin{remark}
Propositions \ref{SemiLagProp1} and  \ref{SemiLagProp2}  are generalization of the results of theorem 5.1 and 6.1 in \cite{debrabant2014semi} in the almost monotone case.
\end{remark}
\begin{remark}
The bound given by the framework is not tight because  we  introduced a switching system in our approach and not in the latter proposition.
\end{remark}
\subsection{Finite Difference scheme}
In this part we suppose that  $a^{\alpha}$ is independent of $x$.
The matrix $a^{\alpha}$ can be written $ P^{\alpha}D^{\alpha}(P^{\alpha})^{t}$ where $P^{\alpha}$ is a unitary matrix with i-th columns $\xi_i^{\alpha}$ and $D^{\alpha}$ is a diagonal $D_{i,j}^{\alpha} = \delta_{i,j} d_i^{\alpha} \ge 0$.
We note $(b^{\alpha})^+(x)$  the vector such that $(b^{\alpha})^+_i(x) = \max( 0, b^{\alpha_i}(x))$ and $(b^{\alpha})^-_i(x) = \max( 0, -b^{\alpha_i}(x))$.
The operator $\mathcal{L}^{\alpha}$ can be discretized for a regular function $u$ using two parameters $h$ and $\hat h$  by 
\begin{eqnarray}
\mathcal{L}^{\alpha}(x, u, \mathcal{D} u, \mathcal{D}^2 u) & \simeq &  \sum_{i=1}^N \frac{d^{\alpha}_i}{2} \frac{2 u(x) - u(x- h \xi_i^{\alpha}) - u(x +h \xi_i^{\alpha})}{h^2}   \nonumber \\
 && - \sum_{i=1}^N (b^{\alpha})^+_i(x) \frac{u(x + \hat h e_i) -u(x)}{\hat h} + \sum_{i=1}^N (b^{\alpha})^{-}_i(x) \frac{u(x ) -u(x- \hat h e_i)}{\hat h}  \nonumber \\
& &  +c^{\alpha}(x) u(x) -f(x),  \nonumber
\end{eqnarray}
where $(e_i)_{i=1,N}$ is the canonical basis in $\R^N$.  We suppose that the equation has been normalized such that
\begin{eqnarray}
\sup_{\alpha} \sum_{i} d_i^{\alpha} + |b^{\alpha}_i(x)|  \le 1. \nonumber
\end{eqnarray}
This is always possible because of \ref{refAss1} and noting that  $\sum_{i} d_i^{\alpha}=tr[a^{\alpha}]$.\\
For $z \in \mathcal{Q}^{\alpha} = \{ 0, \left \{ \xi_i^{\alpha} \right \}_{i=1,N},\left\{e_i \right\}_{i=1,N} \}$ we define the transition probability (see \cite{kushner2013numerical}) for
$\frac{h^2}{\hat h} \le 1$ 
\begin{equation}
\begin{array}{l}
p^{\alpha}(x,x)  =  1 -\sum_{i} \left \{ d_i^{\alpha} + |b^{\alpha}_i(x)|  \frac{h^2}{\hat h} \right \}, \nonumber \\
p^{\alpha}(x,x  \pm \hat h e_i)   =  (b^{\alpha})^{\pm}_i(x) \frac{h^2}{\hat h}, \nonumber \\
p^{\alpha}(x,x \pm h \xi_i^{\alpha})      =  \frac{d_i^{\alpha}}{2}.  \nonumber 
\end{array}
\end{equation}
For any bounded continuous function $\phi$, $x,y,z \in \R^N$, we define the operator $\hat S$
\begin{eqnarray}
\label{scheme0B}
\hat S(h,y,t, \phi_x) =  \sup_{\alpha \in \mathcal{A}} \left \{ -\frac{1}{h^2}  \left[ \sum_{z \in \mathcal{Q}^{\alpha}}  p^\alpha(y,y+z) \phi(x+z) -t \right] + c^{\alpha}(y) t - f^{\alpha}(y) \right \}.  \nonumber
\end{eqnarray}
First we are interested in getting a  solution $v_h$ of 
\begin{eqnarray}
\hat S(h,x,v_h(x), (v_h)_x) = 0, x \in \R^N,
\label{SemiDisDF}
\end{eqnarray} 
in $C^{0,1}( \R^N)$ and get a convergence bound of $|u - v_h|$.\\
We introduce another assumption stronger than assumption \ref{refAss2}:
\begin{Aassumption}
\label{refAss3} 
The constant $\lambda$ in \ref{refAss1} satisfies $\lambda  > \max \left (\sup_{\alpha}{ \frac{1}{2} [\sigma^{\alpha}]^2_1 + [b^{\alpha}]_1},  2 \sqrt{N} \sup_{\alpha} {[b^{\alpha}]_1} \right).$ 
\end{Aassumption}
By using \cite{barles2002convergence} we get the  following proposition
\begin{proposition}
Under assumptions  \ref{refAss1}, \ref{refAss3} the  operator $\hat S$ is monotone, consistent, and there exists a unique solution $v_h$ bounded uniformly in $C^{0,1}(\R^N)$ of 
equation \reff{SemiDisDF}
and $u$ solution of \reff{eqHJB} satisfies
\begin{eqnarray}
|u - v_h| & \le & C h^{\frac{1}{2}} \label{EDPCont}.
\end{eqnarray}
\label{finDiffRes}
\end{proposition}
\begin{remark}
We imposed that the $a^{\alpha}$ is independent of $x$ to get a solution $v_h$ in $C^{0,1}$ which is necessary for us in order to satisfy \ref{epsmonotone}.
Without this assumption we can't prove that   $v_h$ is in $C^{0,1}$: only  $L_{\infty}$ bounds on the error are available (see \cite{barles2005error}). 
As an alternative it would have been possible to use Bonnans and al. \cite{bonnans2004fast}   or Krylov \cite{krylov2005rate}  discretization  with Krylov \cite{krylov2005rate} convergence results with 
bounds for $v_h$ in $C^{0,1}(\R^N)$.
\end{remark}
We introduce the operator
\begin{eqnarray}
\label{scheme2B}
S(h,y,t, [\phi]_x) =  \hat S(h, y, t, (\hat I_{\Delta x ,M} \phi)_x).
\end{eqnarray}
We give the convergence results obtained with the Finite Difference scheme : 
\begin{proposition}
Assume \ref{refAss1} hold and that $\Delta x= h^{q}$, $\hat h = h$. 
Then the scheme \reff{scheme2B} satisfies assumptions \ref{epsmonotone}, \ref{regular}, \ref{consistency} with
$k_2 = min(1,2q-2)$, $k_4= 2$ , $p=q-2$ \\
\end{proposition}
\proof
    First assumption \ref{regular} follows easily from \ref{refAss1}. 
\ref{consistency} follows easily using  lemma \ref{interpolerror}: for $v$ regular  the consistency error is bounded  by
\begin{eqnarray}
\epsilon(h,\Delta x) & = C ( h^2 |\mathcal{D}^4 v |_0 + \hat h |\mathcal{D}^2 v |_0 +  \Delta x^2 (\frac{1}{h^2}+\frac{1}{\hat h }) |\mathcal{D}^2 v |_0).
\end{eqnarray}
Using the fact that $\sum_{z \in \mathcal{Q}^{\alpha}}  p^\alpha(y,y+z) =1$  with positive weights,
\begin{eqnarray}
S(h,y,t+m, [ I_{\Delta x ,M}(\phi+m)]_x) & =    &  S(h,y,t, [ I_{\Delta x ,M}(\phi)]_x), \nonumber
\end{eqnarray}
and property \reff{monoton2} is checked.\\
Suppose $v\in C_b(\R^N)$ , $w \in C^{0,1}(\R^N)$,  $v \ge w$, using estimate \reff{estimEsp1}, the  fact that $\sum_{z \in \mathcal{Q}^{\alpha}}  p^\alpha(y,y+z) =1$  with positive weights,
\begin{eqnarray}
S(h,y,t,[ \hat I_{\Delta x ,M}(v)]_x) \le   S(h,y,t,[ \hat I_{\Delta x ,M}(w)]_x)  +   |w|_1 \frac{ \sqrt{N}\Delta x}{h^2}.  \nonumber
\end{eqnarray}
Similarly  if $v\in C^{0,1}(\R^N)$ , $w \in C_b(\R^N)$,  $v \ge w$ using estimate \reff{BestimEsp2}, 
\begin{eqnarray}
S(h,y,t,[ \hat I_{\Delta x ,M}(v)]_x) \le   S(h,y,t,[ \hat I_{\Delta x ,M}(w)]_x)  +     |v|_1 \frac{ \sqrt{N} \Delta x}{h^2},  \nonumber
\end{eqnarray}
so that  the \reff{monoton1} and \reff{monoton1bis} properties are checked.
\qedsymbol\\  

We then have to check that we can find an $\epsilon$ solution $u_h$ satisfying \reff{eqHJBDis}.
We introduce the operator defined for  $U$ a function on $X_{\Delta x, M}$ : for $ x \in X_{\Delta x, M}$
\begin{eqnarray}
(T_{h,\Delta x} U)(x) & =&  \inf_{\alpha \in \mathcal{A}} \left \{ \frac{1}{ 1+ h^2 c^{\alpha}(x)} \sum_{z \in \mathcal{Q}}  p^\alpha(x,x+z)  \hat I_{\Delta x ,M}(U) (x+z)  + h^2 f^{\alpha}(x) \right \} ,  \nonumber
\end{eqnarray}
and $T_{h,\Delta x}^s$ operator is still defined by equation \reff{defineTs}.
\begin{proposition}
\label{epsSolution2}
Assume \ref{refAss1}, \ref{refAss3} hold. Suppose that   $\Delta x= h^{q}$ with $q>2$, $\hat h =h$.
There exists $s \in \N$ depending on $h$ and  $C$ independent of $h$ such that  $u_h = \hat I_{h,\Delta x} (T_{h,\Delta x}^s 0)$ is an $\epsilon(C h^{q-4})$ solution of scheme $S$.
\end{proposition}
\proof
The proof is similar to the one of proposition  \ref{epsSolution1}.
We first prove that
\begin{eqnarray}
|(T_{h,\Delta x}^k 0)(x) -v_h(x)| & \le & \frac{1}{(1+ \lambda h^2)^k} |v_h|_0 + C \frac{\Delta x}{h^2}.
\end{eqnarray}
Taking $k= \min( i \in \N \mbox{ such that } i \ge -\frac{(q-2) log(h)}{log(1+ \lambda h^2)})$ , using the fact that $\Delta x = h^{q}$, we get that
\begin{eqnarray}
|(T_{h,\Delta x}^k 0)(x) -v_h(x)|   \le C h^{q-2}. \label{iterationT2}
\end{eqnarray}
Let's prove that   $u_h = \hat I_{\Delta x,N} (T_{h,\Delta x}^k 0)$ is an $\epsilon(h^{q-4})$ solution of the scheme $S$.\\
As in  \reff{interpest},
\begin{eqnarray}
|\hat I_{\Delta x ,M} (T_{h,\Delta x}^k 0)(x) -v_h(x)| & \le & |T_{h,\Delta x}^k 0- v_h|_0 + C  \Delta x \nonumber \\
& \le &  C h^{q-2}, \nonumber 
\end{eqnarray}
where we have used  \reff{iterationT2}, so
\begin{eqnarray}
|u_h(x) -v_h(x)| \le   C h^{q-2}. \label{directBound2}
\end{eqnarray}
Then using \ref{refAss1}, \reff{directBound2}, the relation $|\sup A - \sup B| \le \sup | A - B|$,
 the fact that the probabilities are between $0$ and $1$ with sum equal to $1$ and the fact that $u_h =  \hat I_{h,\Delta x} u_h$ :
\begin{eqnarray}
| S(h,x,u_h(x), [u_h]_x)|  & =  & | \hat S(h,x,u_h(x), (u_h)_x) - \hat S(h,x,v_h(x), (v_h)_x)|  \nonumber \\
                          & \le  &    \frac{1}{h^2} |u_h -v_h|_0   \nonumber \\
                          & \le & c h^{q-4}.
\end{eqnarray}

\qedsymbol \\
\begin{proposition}
Assume \ref{refAss1}, \ref{refAss3}, \ref{discommand}  hold.  Constructing an $\epsilon$ solution of \reff{eqHJBDis} with the $u_h$ given by proposition \ref{epsSolution2},
we get 
\begin{eqnarray}
|u- u_h|_0 & \le &  h^{\frac{1}{5}}.
\end{eqnarray}
with $q$ above $21/5$.
\end{proposition}
\proof
The rate of convergence and the value $q=\frac{21}{5}$ is due to a direct use of theorem \ref{mainTheo}.
\qedsymbol \\

\begin{proposition}
Assume \ref{refAss1}, \ref{refAss3} hold. Using  $u_h$ given by proposition \ref{epsSolution2}, we get
\begin{eqnarray}
|u- u_h|_0 & \le &  h^{\frac{1}{2}}.
\end{eqnarray}
\end{proposition}
\proof
This is a  direct use of \reff{directBound2} with the rate of convergence $|u-v_h|_0 \le h^{\frac{1}{2}}$  of proposition \ref{finDiffRes} taking $q=\frac{5}{2}$.
\qedsymbol\\

\section{Numerical tests}
The theoretical bounds obtained in the previous section are not better than the ones obtained with linear interpolation.
The interest of this approximation relies on the fact that where the solution is smooth we expect that the solution won't be truncated and that the
consistency error will be far better than the theoretical one.
All case treated are two dimensional  cases.  For the first three examples, we only give results for the Semi-Lagrangian scheme because the Finite Difference scheme
developed coincides with the classical Finite Difference.
The domain linked to the resolution of equation  \reff{eqHJB} will be noted $Q$, $\mathbb{1}$ is the diagonal unitary matrix and $\mathop{1}$ is the vector
with $1$ components. 
For all Semi-Lagrangian Schemes we choose $h=0.0002$. We discretized the one dimensional space of controls $\mathcal{A}$ with  $2000$ controls.
The software is parallelized with 48 cores as explained in \cite{warin2016some}.
The interpolation used are either linear (2 points per mesh in each direction, monotone scheme), or quadratic (3 points per mesh in each direction), or Cubic (4 points
per mesh in each direction). 
On all the cases and all the tests the fixed point iteration is converging but it can be very slow especially for Finite Differences.
The maximum number of iterations is taken equal to $100000$ and no acceleration was used. The stopping criterion  between
iteration $i$ and $i+1$ was  $$|u^{i+1} -u^i| \le 10^{-7}.$$ The different schemes  are stable and numerically convergent. 
The numerical order of convergence does not show a regular behavior so it
has not been reported.
In the table Err is for the $\mathcal{L}^\infty$ norm for a given discretization given by a number of mesh ( the same in each direction), while  ItN gives the number 
of fixed point iteration used and the computational time is given in seconds.
\subsection{A first regular test case}
The solution of this test case is given by 
\begin{eqnarray}
u(x,y) &=& \sin( \pi x) \sin(\pi y). \nonumber 
\end{eqnarray}
The coefficients are given by
\begin{eqnarray}
  c_a(t,x)  & = &  C , \quad \sigma_a(t,x) = \sigma a  \mathbb{1}   ,  \quad b^{\alpha} = b \mathop{1}, \nonumber
\end{eqnarray}
and the function $f^{\alpha}$ is given by 
\begin{eqnarray}
f^{\alpha}(x,y) & =&  (C + \pi^2 \sigma^2 1_{u(x,y)>0}) u(x,u) - b \pi(\cos(\pi x) \sin(\pi y) + \sin(\pi x) \cos(\pi y)). \nonumber 
\end{eqnarray}
The numerical coefficients are 
\begin{eqnarray}
Q =   [0,2]^2, \quad  b= 0.3, \quad  c = 0.55, \quad \sigma =1 ,  \quad \mathop{A} =[ 0,1],   \nonumber
\end{eqnarray}
and we use a zero Dirichlet boundary condition.
Results are given in table \ref{case1} and this first regular case clearly indicates that the quadratic approximation is by far the
most efficient interpolation : even on this regular case the use of a cubic interpolator doesn't decrease the error with   a CPU time multiplied  at least three-fold.
 In fact it is rapidly converging to the $h$ discretized operator so that the interpolation
error becomes negligible  for a number of mesh equal to 80.
\begin{table}[h]
\caption{Test case 1} \label{case1}
\centering
\begin{tabular}{|c|c|c|c|c|c|c|c|c|c|}
\hline
 \# mesh & \multicolumn{3}{c|}{LINEAR} & \multicolumn{3}{c|}{QUADRATIC} & \multicolumn{3}{c|}{CUBIC}\\
\hline
     &  Err & ItN & Time &   Err & ItN & Time     &   Err & ItN & Time             \\
10   & 1.169 & 25900 & 138     &  0.051  & 29483 & 919   &   0.183   &  45108 & 5182 \\
20   & 1.028 & 26555 & 468     &  0.0065 & 29398 & 3619  &    0.0083  & 29486 &  10770 \\
40   & 0.758 & 27211 & 1879    & 0.0012 & 29485 & 12633  &   0.0011 &   29706 &  43371  \\
80   & 0.243 & 28076 & 6923    & 0.0003 & 29621 & 50988  &  0.0003  & 29788  & 175163    \\
160  & 0.103 & 28620 & 27762    & 0.0003  & 29748 &19517 &    &  &     \\
320  & 0.018 & 29282 &  108000   &      &     &      &     &   &          \\   \hline
\end{tabular}
 \end{table}
\subsection{A second regular problem}
\label{secondRegular}
The solution is here again
\begin{eqnarray}
u(x,y) &=& \sin( \pi x) \sin(\pi y), \nonumber 
\end{eqnarray}
\begin{eqnarray}
  c_a(t,x)  & = &  C , \quad \sigma_a(t,x) = \sigma a  \mathbb{1}   ,  \quad b^{\alpha} = b ( a, \sqrt{1-a^2}), a \in [ \underbar a, \bar a].  \nonumber
\end{eqnarray}
Noting $$\tilde a = \frac{\sin(\pi y) \cos(\pi x)}{\sqrt{ \sin(\pi x)^2 \cos(\pi y)^2 + \cos(\pi x)^2 \sin(\pi y)^2}},$$
$$\phi(a) = a \sin(\pi y) \cos(\pi x) + \sqrt{1- a^2} \sin(\pi x) \cos(\pi y),$$ for $b \le 0,$
the function $f^{\alpha}$ is here given by
\begin{eqnarray}
f^{\alpha}(x,y) &=&  (C + \pi^2 \sigma^2 ) u(x,u) -b \pi K,  \nonumber 
\end{eqnarray}
where $K$ is the maximum of $\phi(\underbar a)$, $\phi(\bar a)$ and $\phi(\tilde a)$ conditionally to  $ \underbar a \le \tilde a \le \bar a$.
We take $ \underbar a =-1$, $\bar a=1$, $\sigma = 1$, $C = 0.6$, $b=-1$, $Q=[0,\frac{1}{2}]$. As boundary condition we take the Dirichlet value given by $u$.\\
Results obtained in table \ref{case2} still show the superiority of the quadratic interpolation for this regular problem with exactly the same conclusions.

\begin{table}[h]
\caption{Test case 2} \label{case2}
\centering
\begin{tabular}{|c|c|c|c|c|c|c|c|c|c|}
\hline
 \# mesh & \multicolumn{3}{c|}{LINEAR} & \multicolumn{3}{c|}{QUADRATIC} & \multicolumn{3}{c|}{CUBIC}\\
\hline
     &  Err & ItN & Time &   Err & ItN & Time & Err & ItN & Time \\
8   & 0.1195 & 493 & 2 &  0.0115  & 1304 & 33      &  0.012 &  1571 &   99       \\
16   & 0.0632 & 903 & 13 & 0.0022  & 1318  & 100   &  0.0022 &  1589&  390          \\
32   & 0.0191 & 1197 & 50 &  0.0003 & 1319 & 405   &  0.0003 &  1590 & 1574            \\
64   & 0.0062 & 1278 & 207 &       &       &        &        &       &            \\
\hline
\end{tabular}
 \end{table}

\subsection{A non regular problem}
We keep the same notations as in the subsection \ref{secondRegular}.\\
Introducing
$$\tilde a = \frac{\frac{1}{2}\sin(\pi y) \cos(\frac{1}{2} \pi x)}{\sqrt{ \sin( 0.5 \pi x)^2 \cos(\pi y)^2 + \frac{1}{4}\cos( 0.5 \pi x)^2 \sin(\pi y)^2}},$$
$$ \hat a = -\tilde a,$$
$$ \hat \phi(a) = \frac{1}{2} a \sin( \pi y) \cos(\frac{1}{2} \pi x) + \sqrt{1- a^2} \sin(\frac{1}{2} \pi x) \cos(\pi y),$$ for $b \le 0,$
the function $f$ is then
\begin{eqnarray}
f^{\alpha} =  \left \{ \begin{array}{ll}
 (C + \pi^2 \sigma^2 ) \sin(\pi y ) \sin( \pi x)  -b \pi K &  \mbox{ for } -1 \le  x \le 0,   \\
 (C +  \pi^2 \sigma^2 \frac{5}{8})  \sin(\frac{1}{2} \pi x) \cos( \pi x)-b \pi \hat K &   \mbox{ for } 0 \le x \le 1,
 \end{array} \right.  \nonumber 
\end{eqnarray}
where $\hat K$ is the maximum between $\hat \phi(\underbar a)$, $\hat \phi(\bar a)$  and  $\hat \phi(\tilde a)$ conditionally to $\underbar a \le \tilde a \le \bar a$, $\hat \phi(\hat a)$ conditionally to $\underbar a \le \hat a \le \bar a$.
We take the values $ \underbar a =-1$, $\bar a=1$, $\sigma = 1$, $C = 0.6$, $b=-1$, $Q=[-1,1]^2$.\\
The boundary conditions are given by the values of the following function  $\psi$:
\begin{eqnarray}
\psi(x,y) =  \sin(\pi y ) \left \{ \begin{array}{ll}
                              \sin( \pi x) & \mbox{  pour }   -1 \le  x \le 0, \\
                              \sin( \frac{1}{2}\pi x) & \mbox{  pour } 1 \ge  x \ge 0.
     \end{array} \right.  \nonumber 
\end{eqnarray}
This test case is interesting  because the continuous function $\psi$ is a regular solution for $x < 0$ and for $x >0$ but it turns out that
it is not the viscosity solution $u$ of the problem.  The reference solution (an estimation of $u$) is numerically calculated with a quadratic interpolation with 128 meshes per directions.
All the methods seems to converge towards the same solution but here again the superiority of the quadratic interpolation is obvious.

\begin{table}[h]
\caption{Test case 3} \label{case3}
\centering
\begin{tabular}{|c|c|c|c|c|c|c|c|c|c|}
\hline
 \# mesh  & \multicolumn{3}{c|}{LINEAR} & \multicolumn{3}{c|}{QUADRATIC} & \multicolumn{3}{c|}{CUBIC} \\
\hline
     &  Err & ItN & Time &   Err & ItN & Time  &   Err & ItN & Time \\
8   & 0.970       & 1810  & 6           &  0.10     & 9406 & 239      & 0.097  & 8877& 536 \\
16   & 0.909      & 3256  & 44          &  0.0141  &  10546 & 777     & 0.0132 & 10677 & 2553 \\
32   & 0.767      &  5551 & 214         &  0.00363 & 10942 & 3277     & 0.0030 & 10916 & 10469  \\
64   & 0.436      & 8259   & 1313       &  0.00187 &  10984   & 12107 & 0.0004 & 10980 & 41952 \\
128  & 0.140      &  8854    & 5348     &        &            &       &        &       &   \\
256  &  0.0478    &  10386    & 26056   &        &            &       &        &       &   \\
\hline
\end{tabular}
 \end{table}

\subsection{A last problem for degeneracy of the diffusion operator}
The last test case will allow us to test the Finite Difference method with a diffusion operator which is degenerated :
The solution is here again
\begin{eqnarray}
u(x,y) &=& \sin( \pi x) \sin(\pi y). \nonumber 
\end{eqnarray}
The coefficients are given by
\begin{eqnarray}
  c_a(t,x)  & = &  C , \quad \sigma_a(t,x) =  \sigma \left( \begin{array}{l} 1 \\ a \end{array} \right)     ,  \quad b^{\alpha} = b \left( \begin{array}{l} 1 \\ a \end{array} \right), \nonumber
\end{eqnarray}
and the function $f^{\alpha}$ is given by 
\begin{eqnarray}
f^{\alpha}(x,y) & =&  \sup_{\hat a \in [ \underline a, \bar a]} \left ((C + \pi^2 \sigma^2 (1+ \hat a^2)) u(x,u)  \right.  \\
& &  \left .- (\sigma^2 \pi^2 (\cos(\pi x) \cos(\pi y) + \pi b \sin(\pi x) \cos(\pi y))) \hat a  \right)\nonumber \\
              &  & -b \pi \cos(\pi x) \sin(\pi y). \nonumber 
\end{eqnarray}
The values taken are  $ \underbar a =-1$, $\bar a=1$, $\sigma = 1$, $C = 0.7$, $b=0.5$, $Q=[-1,1]^2$.
The results clearly indicate that the Finite Difference method proposed is not competitive with  the Semi Lagrangian scheme.
The convergence  of the fixed point iteration for a step $h=0.01$ is not achieved with Finite Difference with $100000$ iterations (the error given between the last two iterations is around $2. 10^{-7}$): there is certainly no exact solution to the scheme.
Once again, the Semi Lagrangian scheme  with quadratic interpolation is the most effective method.\\
\begin{table}[h]
\caption{Test case 4} \label{case4}
\centering
\begin{tabular}{|c|c|c|c|c|c|c|c|c|c|}
\hline
\multicolumn{10}{|c|}{Semi Lagrangian} \\  \hline
 \# mesh  & \multicolumn{3}{c|}{LINEAR} & \multicolumn{3}{c|}{QUADRATIC}  & \multicolumn{3}{c|}{CUBIC} \\
\hline
     &  Err & ItN & Time &   Err & ItN & Time      &      Err & ItN & Time   \\
8   & 0.958 & 1298 & 2 & 0.126   & 12879  & 176    &    0.117 &  17817 &  564                 \\
16   & 0.838 & 3448 & 20 & 0.0227  & 14066 & 576   &   0.0147 &  18854 &   2347              \\
32   & 0.639 & 6841 & 653 & 0.0024  & 14326& 2314  &   0.0015 &   19088&  9430            \\
64  & 0.154 & 13004& 4825 & 0.0018 & 14429&  8499  &          &    &            \\
\hline
\end{tabular}

\begin{tabular}{|c|c|c|c|c|c|c|c|c|c|}
\hline
\multicolumn{10}{|c|}{Finite Difference} \\ \hline
& \multicolumn{9}{c|}{$h= 0.01$}  \\
\hline
 \# mesh & \multicolumn{3}{c|}{LINEAR} & \multicolumn{3}{c|}{QUADRATIC}& \multicolumn{3}{c|}{CUBIC} \\
\hline
     &  Err & ItN & Time &   Err & ItN & Time           &      Err & ItN & Time         \\
8   &  0.9598 & 8258 & 38 &  0.200   & $10^5$  & 1649   &  0.139     &  $10^5$ & 7168         \\
16   & 0.9276 &  15240 & 257 & 0.064  & $10^5$ & 4986   &  0.0345    &  $10^5$ & 28590        \\
32   & 0.8517 &  26815 & 1307 & 0.023  & $10^5$ & 19450 &  0.0129    &  $10^5$ & 113792       \\
64  &  0.69840 & 46830  & 8787 & 0.008      & $10^5$     &  145509      &             &         &     \\
128 &  0.4599 & 79658  & 49851&       &      &     &             &      &   \\
256 &  0.119 & 50000 &  198970 &      &      &     &            &       & \\ \hline
\end{tabular}

\end{table}

\thanks{
Special thank to  the two anonymous referees for their corrections and suggestions.}
\bibliographystyle{alpha}
\bibliography{MyBib}

\newcommand{\etalchar}[1]{$^{#1}$}
\begin{thebibliography}{BFF{\etalchar{+}}15}

\bibitem[AA00]{augoula2000high}
Steeve Augoula and R{\'e}mi Abgrall.
\newblock High order numerical discretization for hamilton--jacobi equations on
  triangular meshes.
\newblock {\em Journal of Scientific Computing}, 15(2):197--229, 2000.

\bibitem[Abg09]{abgrall2009construction}
Remi Abgrall.
\newblock Construction of simple, stable, and convergent high order schemes for
  steady first order hamilton--jacobi equations.
\newblock {\em SIAM Journal on Scientific Computing}, 31(4):2419--2446, 2009.

\bibitem[ADM93]{azaiez1993methodes}
Mejdi Aza{\"\i}ez, Monique Dauge, and Yvon Maday.
\newblock M{\'e}thodes spectrales et des {\'e}l{\'e}ments spectraux.
\newblock 1993.

\bibitem[Atk08]{atkinson2008introduction}
Kendall~E Atkinson.
\newblock {\em An introduction to numerical analysis}.
\newblock John Wiley \& Sons, 2008.

\bibitem[BFF{\etalchar{+}}15]{bokanowski2015value}
Olivier Bokanowski, Maurizio Falcone, Roberto Ferretti, Lars Gr{\"u}ne, Dante
  Kalise, and Hasnaa Zidani.
\newblock Value iteration convergence of $\epsilon$-monotone schemes for
  stationary hamilton-jacobi equations.
\newblock {\em Discrete and Continuous Dynamical Systems-Series A},
  35(9):4041--4070, 2015.

\bibitem[BFS16]{bokanowski2016efficient}
Olivier Bokanowski, Maurizio Falcone, and Smita Sahu.
\newblock An efficient filtered scheme for some first order time-dependent
  hamilton--jacobi equations.
\newblock {\em SIAM Journal on Scientific Computing}, 38(1):A171--A195, 2016.

\bibitem[BJ02]{barles2002convergence}
Guy Barles and Espen~Robstad Jakobsen.
\newblock On the convergence rate of approximation schemes for
  hamilton-jacobi-bellman equations.
\newblock {\em ESAIM: Mathematical Modelling and Numerical Analysis},
  36(1):33--54, 2002.

\bibitem[BJ05]{barles2005error}
Guy Barles and Espen~R Jakobsen.
\newblock Error bounds for monotone approximation schemes for
  hamilton--jacobi--bellman equations.
\newblock {\em SIAM journal on numerical analysis}, 43(2):540--558, 2005.

\bibitem[BMZ10]{bokanowski2010convergence}
Olivier Bokanowski, Nadia Megdich, and Hasnaa Zidani.
\newblock Convergence of a non-monotone scheme for hamilton--jacobi--bellman
  equations with discontinous initial data.
\newblock {\em Numerische Mathematik}, 115(1):1--44, 2010.

\bibitem[BOZ04]{bonnans2004fast}
J~Fr{\'e}d{\'e}ric Bonnans, {\'E}lisabeth Ottenwaelter, and Housnaa Zidani.
\newblock A fast algorithm for the two dimensional hjb equation of stochastic
  control.
\newblock {\em ESAIM: Mathematical Modelling and Numerical Analysis},
  38(4):723--735, 2004.

\bibitem[BPR16]{bokanowski2016high}
Olivier Bokanowski, Athena Picarelli, and Christoph Reisinger.
\newblock High-order filtered schemes for time-dependent second order hjb
  equations.
\newblock {\em arXiv preprint arXiv:1611.04939}, 2016.

\bibitem[BS91]{barles1991convergence}
Guy Barles and Panagiotis~E Souganidis.
\newblock Convergence of approximation schemes for fully nonlinear second order
  equations.
\newblock {\em Asymptotic analysis}, 4(3):271--283, 1991.

\bibitem[CF95]{camilli1995approximation}
Fabio Camilli and Maurizio Falcone.
\newblock An approximation scheme for the optimal control of diffusion
  processes.
\newblock {\em ESAIM: Mathematical Modelling and Numerical Analysis},
  29(1):97--122, 1995.

\bibitem[DJ13]{debrabant2013semi}
Kristian Debrabant and Espen Jakobsen.
\newblock Semi-lagrangian schemes for linear and fully non-linear diffusion
  equations.
\newblock {\em Mathematics of Computation}, 82(283):1433--1462, 2013.

\bibitem[DJ14]{debrabant2014semi}
Kristian Debrabant and Espen~R Jakobsen.
\newblock Semi-lagrangian schemes for linear and fully non-linear
  hamilton-jacobi-bellman equations.
\newblock {\em arXiv preprint arXiv:1403.1217}, 2014.

\bibitem[EF79]{evans1979optimal}
Lawrence~C Evans and Avner Friedman.
\newblock Optimal stochastic switching and the dirichlet problem for the
  bellman equation.
\newblock {\em Transactions of the American Mathematical Society},
  253:365--389, 1979.

\bibitem[FF13]{falcone2013semi}
Maurizio Falcone and Roberto Ferretti.
\newblock {\em Semi-Lagrangian approximation schemes for linear and
  Hamilton-Jacobi equations}, volume 133.
\newblock SIAM, 2013.

\bibitem[FO13]{froese2013convergent}
Brittany~D Froese and Adam~M Oberman.
\newblock Convergent filtered schemes for the monge--ampère partial
  differential equation.
\newblock {\em SIAM Journal on Numerical Analysis}, 51(1):423--444, 2013.

\bibitem[JP00]{jiang2000weighted}
Guang-Shan Jiang and Danping Peng.
\newblock Weighted eno schemes for hamilton--jacobi equations.
\newblock {\em SIAM Journal on Scientific computing}, 21(6):2126--2143, 2000.

\bibitem[JS13]{jensen2013convergence}
Max Jensen and Iain Smears.
\newblock On the convergence of finite element methods for
  hamilton--jacobi--bellman equations.
\newblock {\em SIAM Journal on Numerical Analysis}, 51(1):137--162, 2013.

\bibitem[KD13]{kushner2013numerical}
Harold Kushner and Paul~G Dupuis.
\newblock {\em Numerical methods for stochastic control problems in continuous
  time}, volume~24.
\newblock Springer Science \& Business Media, 2013.

\bibitem[Kry00]{krylov2000rate}
NV~Krylov.
\newblock On the rate of convergence of finite-difference approximations for
  bellmans equations with variable coefficients.
\newblock {\em Probability theory and related fields}, 117(1):1--16, 2000.

\bibitem[Kry05]{krylov2005rate}
Nicolai~V Krylov.
\newblock The rate of convergence of finite-difference approximations for
  bellman equations with lipschitz coefficients.
\newblock {\em Applied Mathematics \& Optimization}, 52(3):365--399, 2005.

\bibitem[Lep00]{lepsky2000spectral}
Olga Lepsky.
\newblock Spectral viscosity approximations to hamilton--jacobi solutions.
\newblock {\em SIAM journal on numerical analysis}, 38(5):1439--1453, 2000.

\bibitem[LS95]{lions1995convergence}
P-L Lions and PE~Souganidis.
\newblock Convergence of muscl and filtered schemes for scalar conservation
  laws and hamilton-jacobi equations.
\newblock {\em Numerische Mathematik}, 69(4):441--470, 1995.

\bibitem[Men89]{menaldi1989some}
Jos{\'e}-Luis Menaldi.
\newblock Some estimates for finite difference approximations.
\newblock {\em SIAM journal on control and optimization}, 27(3):579--607, 1989.

\bibitem[Obe06]{oberman2006convergent}
Adam~M Oberman.
\newblock Convergent difference schemes for degenerate elliptic and parabolic
  equations: Hamilton--jacobi equations and free boundary problems.
\newblock {\em SIAM Journal on Numerical Analysis}, 44(2):879--895, 2006.

\bibitem[OS91]{osher1991high}
Stanley Osher and Chi-Wang Shu.
\newblock High-order essentially nonoscillatory schemes for hamilton--jacobi
  equations.
\newblock {\em SIAM Journal on numerical analysis}, 28(4):907--922, 1991.

\bibitem[OS15]{oberman2015filtered}
Adam~M Oberman and Tiago Salvador.
\newblock Filtered schemes for hamilton--jacobi equations: A simple
  construction of convergent accurate difference schemes.
\newblock {\em Journal of Computational Physics}, 284:367--388, 2015.

\bibitem[PFV03]{pooley2003numerical}
David~M Pooley, Peter~A Forsyth, and Ken~R Vetzal.
\newblock Numerical convergence properties of option pricing pdes with
  uncertain volatility.
\newblock {\em IMA Journal of Numerical Analysis}, 23(2):241--267, 2003.

\bibitem[QSS10]{quarteroni2010numerical}
Alfio Quarteroni, Riccardo Sacco, and Fausto Saleri.
\newblock {\em Numerical mathematics}, volume~37.
\newblock Springer Science \& Business Media, 2010.

\bibitem[Shu07]{shu2007high}
Chi-Wang Shu.
\newblock High order numerical methods for time dependent hamilton-jacobi
  equations.
\newblock In {\em Mathematics and computation in imaging science and
  information processing}, pages 47--91. World Scientific, 2007.

\bibitem[SS16]{smears2016discontinuous}
Iain Smears and Endre S{\"u}li.
\newblock Discontinuous galerkin finite element methods for time-dependent
  hamilton--jacobi--bellman equations with cordes coefficients.
\newblock {\em Numerische Mathematik}, 133(1):141--176, 2016.

\bibitem[War16]{warin2016some}
Xavier Warin.
\newblock Some non-monotone schemes for time dependent
  hamilton--jacobi--bellman equations in stochastic control.
\newblock {\em Journal of Scientific Computing}, 66(3):1122--1147, 2016.

\end{thebibliography}
\end{document}